\documentclass[12pt]{amsart}

\usepackage{amsmath}
\usepackage{amssymb}
\usepackage{latexsym}
\usepackage{euscript}
\usepackage{eufrak}
\usepackage{verbatim,graphicx}
\usepackage{subfigure}
\usepackage{curves}
\usepackage{bigints}

\newtheorem{theorem}{Theorem}[section]

\newtheorem{corollary}[theorem]{Corollary}

\numberwithin{equation}{section}
\numberwithin{figure}{section}

\theoremstyle{definition}

\newtheorem{definition}[theorem]{Definition} 
\newtheorem{example}[theorem]{Example}

\newtheorem{exercise*}[theorem]{*Exercise}
\newtheorem{remark}[theorem]{Remark}

\theoremstyle{remark}

\def\&{\wedge}
\newcommand{\w}{\omega}
\newcommand{\R}{\mathbb{R}}
\newcommand{\A}{\mathbb{A}}
\newcommand{\E}{\mathbb{E}}
\newcommand{\e}{\mathbf{e}}

\newcommand{\calF}{\mathcal{F}}

\newcommand{\bx}{\mathbf{x}}

\begin{document}

\title{A characterization of hyperbolic affine flat, affine minimal surfaces in $\A^3$}

\author{Jeanne N. Clelland}
\address{Department of Mathematics, 395 UCB, University of
Colorado, Boulder, CO 80309-0395}
\email{Jeanne.Clelland@colorado.edu}

\author{Jonah M. Miller}
\address{Department of Mathematics, 395 UCB, University of
Colorado, Boulder, CO 80309-0395}
\email{Jonah.Miller@colorado.edu}

\thanks{This research was supported in part by NSF grant DMS-1206272.}

\begin{abstract}
We investigate the geometric properties of hyperbolic affine flat, affine minimal surfaces in the equiaffine space $\A^3$. We use Cartan's method of moving frames to compute a complete set of local invariants for such surfaces. Using these invariants, we give a complete local classification of such surfaces and construct new examples.  
\end{abstract} 

\maketitle

\section{Introduction}\label{intro-sec}

In equiaffine geometry, one of the most-studied categories of surfaces is the class of {\em affine spheres}.  A nondegenerate surface $\Sigma \subset \A^3$ is called a {\em proper affine sphere} if the affine normal lines passing through each point of $\Sigma$ intersect in a single point, and an {\em improper affine sphere} if the affine normal lines passing through each point of $\Sigma$ are all parallel.  These surfaces are much more plentiful in equiaffine geometry than in Euclidean geometry, where the proper and improper ``spheres" are simply the spheres and planes, respectively.  

Much of the study of improper affine spheres has been devoted to surfaces in the elliptic category.  Any elliptic improper affine sphere can be represented locally as the graph of a solution to the elliptic Monge-Amp\`ere equation
\[ z_{xx} z_{yy} - (z_{xy})^2 = 1, \]
and such surfaces can be given a conformal representation which is very useful in studying their geometric properties; see, e.g., \cite{ACG07}, \cite{FMM96}, \cite{Ferrer02}.  Hyperbolic improper affine spheres have received comparably little attention, but a few results are known: Gao \cite{Gao93} classified all polynomials whose graphs are improper affine spheres without regard to type, and Magid and Ryan \cite{MR90} gave classifications for both elliptic and hyperbolic improper affine spheres under the additional condition that the affine Gauss curvature vanishes identically.

An improper affine sphere $\Sigma$ necessarily has affine mean curvature
$H_{\text{aff}} = 0$.
In particular, an improper affine sphere is also an {\em affine minimal} surface---another category of surfaces which has been the object of considerable study in affine geometry (see, e.g., \cite{Buyske92}, \cite{CT80}, \cite{Krauter94}, \cite{VV89}). If, in addition, the affine metric of $\Sigma$ has Gauss curvature $K_{\text{aff}} = 0$, then $\Sigma$ is called {\em affine flat}.  
In the elliptic category, an affine flat, affine minimal surface must be an improper affine sphere; in fact, it is shown in \cite{MR90} that such a surface must be contained in a paraboloid.

By contrast, in the hyperbolic category there exist surfaces with $K_{\text{aff}} = H_{\text{aff}} = 0$ which are {\em not} improper affine spheres.  Such surfaces may be of independent interest; for example, they are singled out in \cite{CT80} as a special case of affine minimal surfaces for which the surface transformation described there has a particularly simple form.  In this paper, we will use Cartan's method of equivalence to give a complete local classification of hyperbolic affine flat, affine minimal surfaces.  In the process, we recover the classification of hyperbolic, affine flat improper affine spheres given in \cite{MR90}, depending on one arbitrary function of one variable.  In addition, we find a larger family of hyperbolic affine flat, affine minimal surfaces which are {\em not} improper affine spheres, depending on two arbitrary functions of one variable.

This paper is organized as follows: in \S \ref{background-sec} we review the necessary concepts in equiaffine geometry, including the notion of unimodular frames on the equiaffine space $\A^3$ and their associated Maurer-Cartan forms.  In \S \ref{equivalence-prob-sec} we carry out the equivalence method to compute local invariants for hyperbolic affine flat, affine minimal surfaces in $\A^3$.  In \S \ref{normal-forms-sec} we derive a local normal form for a compatible, overdetermined PDE system whose solutions give rise to parametrized surfaces of this type.  In \S \ref{examples-sec} we use solutions of this system to construct examples of such surfaces, and in \S \ref{conclusion-sec} we make some concluding remarks.

\begin{comment}
Our classification reveals the following interesting features:
\begin{itemize}
\item  As any surface in $\A^3$ may locally be represented as the graph of a single function $z=f(x,y)$, the affine flat, affine minimal condition $K_{\text{aff}}=H_{\text{aff}}=0$ may be regarded as an overdetermined PDE system consisting of two PDE's for the single function $z(x,y)$; as such, we might expect it to have no---or at least very few---solutions. (In the the more familiar case of surfaces in Euclidean space $\E^3$, for instance, the only flat, minimal surfaces are planes.)  However, we find a much larger family of solutions than in the Euclidean case: specifically, we find that the family of hyperbolic affine flat, affine minimal surfaces in $\A^3$ is locally parametrized by two arbitrary functions $f(v), \ell(v)$ of a single variable $v$.  
\item The surface associated to given functions $f(v), \ell(v)$ is closely related to the Sturm-Liouville operator $\left(\frac{d}{dv}\right)^2 - \ell(v)$; specifically, the ODE system whose solutions give rise to surfaces of this type contains a subsystem that is equivalent to the Sturm-Liouville equation $h''(v) - \tfrac{1}{2} \ell(v) h(v) = 0$.
\end{itemize}

\end{comment}

\section{Equiaffine space, unimodular frames, and Maurer-Cartan forms}\label{background-sec}

We begin by recalling the definition of equiaffine space $\A^3$ and its symmetry group $A(3)$.  (For a comprehensive introduction to affine geometry, see, e.g., \cite{NS94}.)

\begin{definition}
Three-dimensional {\em equiaffine space} $\A^3$ (which for convenience we will refer to simply as ``affine space") consists of the vector space $\R^3$, together with a nondegenerate volume form
\[ dV: \Lambda^3 \R^3 \to \R. \]
The {\em equiaffine group} $A(3)$ is the group of all transformations $\phi:\E^3 \to \E^3$ which preserve the volume form; it consists of all transformations of the form
\[ \phi(\bx) = A\bx + \mathbf{b}, \]
where $A \in SL(3)$ and $\mathbf{b} \in \A^3$.
\end{definition}

As a vector space, $\A^3$ is equivalent to the Euclidean space $\E^3$.
But while the inner product structure on $\E^3$ induces a volume form on $\E^3$, the converse is false: there is no inner product on $\A^3$ which is preserved by the action of the equiaffine group $A(3)$.  Thus in equiaffine geometry, there are no obvious analogs of metric notions such as length or angles defined on tangent vectors.

We will use Cartan's method of moving frames to compute local invariants for surfaces in $\A^3$.  The notions of ``hyperbolic" (vs. ``elliptic") surfaces, ``affine Gauss curvature", and ``affine mean curvature" will arise during the frame adaptation process, and we will give precise definitions for these terms as we encounter them.

In Euclidean geometry, one usually considers the set of orthonormal frames $(\e_1, \e_2, \e_3)$ for the tangent space $T_{\bx}\E^3$ based at each point $\bx \in \E^3$.  But in equiaffine geometry, there is no well-defined notion of an angle between tangent vectors, and hence no notion of ``orthonormal." Instead, we consider the set of {\em unimodular} frames.

\begin{definition}
A basis $(\e_1, \e_2, \e_3)$ for the tangent space $T_{\bx}\A^3$ at a point $\bx \in \A^3$ is called a {\em unimodular frame} at $\bx$ if 
\[ dV(\e_1, \e_2, \e_3) = 1. \]
\end{definition}
This is equivalent to the condition that the vectors $(\e_1, \e_2, \e_3)$ span a parallelepiped of (oriented) volume 1, and also to the condition that
\begin{equation}
\e_1 \wedge \e_2 \wedge \e_3  = \underline{\e}_1 \wedge \underline{\e}_2 \wedge \underline{\e}_3, \label{unimodular-wedge-condition}
\end{equation}
where $(\underline{\e}_1, \underline{\e}_2, \underline{\e}_3)$ is the standard (constant) basis on $\R^3$.

The unimodular frames on $\A^3$ form a principal fiber bundle $\pi: \mathcal\calF(\A^3) \to \A^3$, with structure group equal to $SL(3)$, called the {\em unimodular frame bundle} over $\A^3$.  The bundle $\calF(\A^3)$ is isomorphic to the affine group $A(3)$.

The {\em Maurer-Cartan forms} $\w^i, \w^i_j$ on $\calF(\A^3)$ are the 1-forms on $\calF(\A^3)$ defined by the equations
\begin{align}
d\bx & = \e_i \w^i, \label{MC-forms-def}   \\ 
d\e_j & =  \e_i \w^i_j. \notag
\end{align}
(Note that we use the Einstein summation convention, and all indices are summed from 1 to 3.)  The 1-forms $\w^1, \w^2, \w^3$ are called the {\em dual forms} (or sometimes the {\em solder forms}), while the 1-forms $\{\w^i_j, \ 1 \leq i,j \leq 3\}$ are called the {\em connection forms}.  They satisfy the Maurer-Cartan structure equations
\begin{align}
d\w^{i} & = - \w^i_j \& \w^j,  \label{structure-eqns} \\
d\w^i_j & = - \w^i_k \& \w^k_j. \notag
\end{align}
(See \cite{IL03} for a discussion of Maurer-Cartan forms and their structure equations.)  Differentiating the relation \eqref{unimodular-wedge-condition}
 yields the relation
\begin{equation}
\w^1_1 + \w^2_2 + \w^3_3 = 0. \label{MC-forms-trace-relation} 
\end{equation}
Unlike in Euclidean geometry, where $\w^i_j = -\w^j_i$, the connection forms on $\calF(\A^3)$ are linearly independent except for the single relation \eqref{MC-forms-trace-relation}.

\section{Equivalence problem and local invariants}\label{equivalence-prob-sec}

In this section, we use Cartan's method of equivalence to construct adapted frames and compute local invariants for hyperbolic surfaces in $\A^3$; in particular, the affine Gauss and mean curvatures $K_{\text{aff}}, H_{\text{aff}}$ will be introduced.

\subsection{Adapted frames on $\Sigma$ and the 0-adapted frame bundle}

Now let $\Sigma \subset \A^3$ be a regular surface in $\A^3$. 

\begin{definition}
The subset $\calF(\Sigma) \subset \calF(\A^3)$ consisting of all unimodular frames based at all points $\bx \in \Sigma$ will be called the {\em unimodular frame bundle} over $\Sigma$.  (Technically, $\calF(\Sigma)$ is the pullback of $\calF(\A^3)$ to $\Sigma$ via the inclusion map $\iota:\Sigma \to \A^3$.)  A section $\sigma:\Sigma \to \calF(\Sigma)$ is called a {\em unimodular frame field} on $\Sigma$.
\end{definition}

In order to reduce notational clutter, for the remainder of the paper we will abuse notation slightly by using $(\e_1, \e_2, \e_3)$ to denote a unimodular frame field $\sigma(\bx) = (\e_1(\bx), \e_2(\bx), \e_3(\bx))$ on $\calF(\Sigma)$.  It should be clear from context when this notation refers to a frame field on $\Sigma$ rather than to a point of $\calF(\Sigma)$.  Furthermore, we will denote the pullbacks $\sigma^*\w^i, \sigma^*\w^i_j$ of the Maurer-Cartan forms to $\Sigma$ by $\bar{\w}^i, \bar{\w}^i_j$, respectively.  While the Maurer-Cartan forms $\w^i, \w^i_j$ are linearly independent 1-forms on $\calF(\A^3)$ (except for the relation \eqref{MC-forms-trace-relation}), the forms $\bar{\w}^i, \bar{\w}^i_j$ on $\Sigma$ are all sections of the rank 2 cotangent bundle $T^*\Sigma$; hence there are many linear dependence relations among them, and these will become apparent during the frame adaptation process.

The method of equivalence begins by considering those unimodular frame fields on $\Sigma$ which are ``nicely" adapted to the geometry of $\Sigma$.  In Euclidean geometry, one typically considers orthonormal frame fields for which the frame vectors $\e_1, \e_2$ are tangent to $\Sigma$ and $\e_3$ is normal to $\Sigma$.  In equiaffine geometry, we have no obvious notion of a ``normal vector" to $\Sigma$, but the concept of tangency is still well-defined.  Thus we will initially consider the following class of unimodular frames on $\Sigma$:

\begin{definition}
A unimodular frame $(\e_1, \e_2, \e_3)$ based at a point $\bx \in \Sigma$ will be called {\em 0-adapted} if the frame vectors $\e_1, \e_2$ span the tangent space $T_{\bx}\Sigma$.  A unimodular frame field $(\e_1, \e_2, \e_3)$ on $\Sigma$ will be called {\em 0-adapted} if, for each $\bx \in \Sigma$, the frame $(\e_1(\bx), \e_2(\bx), \e_3(\bx))$ is a 0-adapted frame at $\bx$. 
\end{definition}

The 0-adapted frame fields on $\Sigma$ are the sections of a subbundle $\calF_0 \subset \calF(\Sigma)$, called the {\em 0-adapted frame bundle}.  
Any two 0-adapted frames $(\e_1, \e_2, \e_3)$, $(\tilde{\e}_1, \tilde{\e}_2, \tilde{\e}_3)$ based at a point $\bx \in \Sigma$ must have the property that
\[ \text{span}(\tilde{\e}_1, \tilde{\e}_2) = \text{span}(\e_1, \e_2); \]
therefore they must differ by a transformation of the form
\begin{equation}
\begin{bmatrix} \tilde{\e}_1 & \tilde{\e}_2 & \tilde{\e}_3 
\end{bmatrix} = \begin{bmatrix} \e_1 & \e_2 & \e_3 \end{bmatrix} 
\begin{bmatrix} B & \begin{matrix} r_1 \\[0.05in] r_2 \end{matrix} 
\\[0.2in] 
\begin{matrix} 0 \ &\ 0 \ \end{matrix} & (\det B)^{-1} \end{bmatrix}, \label{0-adapted-structure-group}
\end{equation}
where $B \in GL(2)$ and $r_1, r_2 \in \R$. Furthermore, if $(\e_1, \e_2, \e_3)$ is any 0-adapted frame on $\Sigma$, then any frame $(\tilde{\e}_1, \tilde{\e}_2, \tilde{\e}_3)$ given by \eqref{0-adapted-structure-group} is also 0-adapted.
The 0-adapted frame bundle $\calF_0$ is a principal bundle, with structure group  $G_0 \subset SL(3)$ equal to the subgroup of $SL(3)$ consisting of all matrices of the form in equation \eqref{0-adapted-structure-group}.

Now consider the pullbacks of equations \eqref{MC-forms-def} to $\Sigma$ via a section of $\calF_0$.  (More intuitively, this means that we now regard $(\e_1, \e_2, \e_3)$ as a 0-adapted frame field on $\Sigma$ and replace the forms $\w^i, \w^i_j$ in equations \eqref{MC-forms-def} with the forms $\bar{\w}^i, \bar{\w}^i_j$ associated with this frame field on $\Sigma$.) 
In particular, from the equation
\[ d\bx = \e_1 \bar{\w}^1 + \e_2 \bar{\w}^2 + \e_3 \bar{\w}^3 \]
and the fact that the image of $d\bx$ spans the tangent space $T_{\bx}\Sigma$ at each point $\bx \in \Sigma$, it follows that $\bar{\w}^3 = 0$, and that $\bar{\w}^1, \bar{\w}^2$ are linearly independent 1-forms which span the cotangent space $T^*_{\bx}\Sigma$ at each point $\bx \in \Sigma$.

Differentiating the equation $\bar{\w}^3=0$ according to the structure equations \eqref{structure-eqns} implies that
\[ 0 = d\bar{\w}^3 = -(\bar{\w}^3_1 \& \bar{\w}^1 + \bar{\w}^3_2 \& \bar{\w}^2). \]
By Cartan's lemma (see \cite{IL03}), it follows that there exist functions $h_{ij} = h_{ji}$ such that
\begin{equation}
 \begin{bmatrix} \bar{\w}^3_1 \\[0.1in] \bar{\w}^3_2 \end{bmatrix} = \begin{bmatrix} 
h_{11} & h_{12} \\[0.1in] h_{12} & h_{22} \end{bmatrix} \begin{bmatrix} 
\bar{\w}^1 \\[0.1in] \bar{\w}^2 \end{bmatrix}. \label{h-matrix}
\end{equation}

\subsection{Reduction of the structure group}

The method of equivalence proceeds by examining how the functions $h_{ij}$ in equation \eqref{h-matrix} vary among different choices of 0-adapted frame fields on $\Sigma$, and by choosing from among the 0-adapted frames a subset of frames for which the $h_{ij}$ are somehow normalized. Then we look for new relations among the Maurer-Cartan forms associated to this restricted class of adapted frame fields.  This process is then iterated until---hopefully---we arrive at a single, canonical choice of unimodular frame at each point of $\Sigma$.

So suppose that two 0-adapted frame fields $(\e_1, \e_2, \e_3), (\tilde{\e}_1, \tilde{\e}_2, \tilde{\e}_3)$ on $\Sigma$, with associated Maurer-Cartan forms $(\bar{\w}^i, \bar{\w}^i_j), (\tilde{\bar{\w}}^i, \tilde{\bar{\w}}^i_j)$, respectively, are related by a transformation of the form \eqref{0-adapted-structure-group}.  Equations \eqref{MC-forms-def} imply that
\begin{equation}
\begin{bmatrix} \tilde{\bar{\w}}^1 \\[0.1in] \tilde{\bar{\w}}^2 \end{bmatrix} = 
B^{-1} \begin{bmatrix} \bar{\w}^1 \\[0.1in]  \bar{\w}^2 \end{bmatrix}, \qquad \qquad 
\begin{bmatrix} \tilde{\bar{\w}}^3_1 \\[0.1in] \tilde{\bar{\w}}^3_2 \end{bmatrix} = 
(\det B)\, {}^t\hskip-2.5pt B \begin{bmatrix} \bar{\w}^3_1 \\[0.1in] \bar{\w}^3_2 \end{bmatrix} , \label{MC-forms-trans}
\end{equation}
and it follows that the the functions $h_{ij}$ of equation \eqref{h-matrix} and the corresponding functions $\tilde{h}_{ij}$ for the transformed frame field $(\tilde{\e}_1, \tilde{\e}_2, \tilde{\e}_3)$ are related by the equation
\begin{equation}
 \begin{bmatrix} \tilde{h}_{11} & \tilde{h}_{12} \\[0.1in] \tilde{h}_{12} & 
\tilde{h}_{22} \end{bmatrix} = (\det B)\, {}^t\hskip-2.5pt B \begin{bmatrix} 
h_{11} & h_{12} \\[0.1in] h_{12} & h_{22} \end{bmatrix} B  . \label{h-action}
\end{equation}

We may regard equation \eqref{h-action} as defining an action of $GL(2)$ on the space of $2 \times 2$ symmetric matrices $h = [h_{ij}]$.  This action preserves the sign of the determinant of $h$; therefore the sign of $\det(h(\bx))$ is the same for any 0-adapted frame based at a point $\bx \in \Sigma$.  

\begin{definition}
If the matrix $[h_{ij}]$ is nonsingular at every point of a regular surface $\Sigma$, then $\Sigma$ is called {\em nondegenerate}.  Furthermore, a nondegenerate surface $\Sigma$ is called:
\begin{itemize}
\item {\em elliptic} if $\det[h_{ij}] >0$ at every point of $\Sigma$;
\item {\em hyperbolic} if $\det[h_{ij}] <0$ at every point of $\Sigma$.
\end{itemize}
\end{definition}

\begin{remark}
The sign of $\det(h)$ is the same as the sign of the Gauss curvature $K$ of $\Sigma$ when regarded as a surface in Euclidean space $\E^3$.  (While the Gauss curvature of $\Sigma$ is not invariant under the group of equiaffine transformations, its {\em sign} is well-defined up to equiaffine transformations.)  Thus this division of nondegenerate surfaces into elliptic and hyperbolic types is, in fact, equivalent to the usual Euclidean notions of elliptic ($K>0$) and hyperbolic ($K<0$) surfaces. (See \cite{NS94} for details.)
\end{remark}

For the remainder of this paper, we will assume that $\Sigma$ is hyperbolic.

The $h_{ij}$ are real-valued functions on the 0-adapted frame bundle $\calF_0$ of $\Sigma$, and the
$GL(2)$-action \eqref{h-action} is transitive on the set of all $2 \times 2$ symmetric matrices of negative determinant. Therefore, there exists a nonempty subbundle $\calF_1 \subset \calF_0$ consisting of those 0-adapted frames on $\Sigma$ for which
\begin{equation}
\begin{bmatrix} h_{11} & h_{12} \\[0.1in] h_{12} & h_{22} \end{bmatrix} = \begin{bmatrix} 0\ &\ 1 \\[0.1in] 1\ &\ 0 \end{bmatrix}. \label{h-1-adapted}
\end{equation}

\begin{definition}
The bundle $\calF_1$ will be called the {\em 1-adapted frame bundle} on $\Sigma$. 
Any frame $(\e_1, \e_2, \e_3) \in \calF_1$ will be called a {\em 1-adapted frame} on $\Sigma$, and any section of $\calF_1$ will be called a {\em 1-adapted frame field} on $\Sigma$.
\end{definition}

Equations \eqref{h-matrix} and \eqref{h-1-adapted} imply that for any 1-adapted frame field on $\Sigma$, the associated Maurer-Cartan forms satisfy the relations
\begin{equation}
\bar{\w}^3_1 = \bar{\w}^2, \qquad \bar{\w}^3_2 = \bar{\w}^1. \label{1-adapted-MC-relns}
\end{equation}

It is straightforward to show that any two 1-adapted frames $(\e_1, \e_2, \e_3)$, $(\tilde{\e}_1, \tilde{\e}_2, \tilde{\e}_3)$ based at a point $\bx \in \Sigma$ must differ by a transformation of the form 
\begin{equation}
\begin{bmatrix} \tilde{\e}_1 & \tilde{\e}_2 & \tilde{\e}_3 
\end{bmatrix} = \begin{bmatrix} \e_1 & \e_2 & \e_3 \end{bmatrix} 
\begin{bmatrix} \epsilon_1 e^{\lambda} & 0 & r_1 \\[0.1in] 0 & \epsilon_2 e^{-\lambda} & r_2 \\[0.1in] 0 & 0 & \epsilon_1 \epsilon_2 \end{bmatrix}, \label{1-adapted-structure-group}
\end{equation}
where $\lambda, r_1, r_2 \in \R$ and $\epsilon_1, \epsilon_2 = \pm 1$.  The 1-adapted frame bundle $\calF_1$ is a principal bundle, with structure group  $G_1 \subset G_0$ equal to the subgroup of $SL(3)$ consisting of all matrices of the form in equation \eqref{1-adapted-structure-group}.

\begin{definition}
The quadratic form
\begin{equation}
 \text{I}_{\text{aff}} = \bar{\w}^3_1 \circ \bar{\w}^1 + \bar{\w}^3_2 \circ \bar{\w}^2  \label{first-fundamental form}
\end{equation}
on the 1-adapted frame bundle $\calF_1$ is called the {\em affine first fundamental form} of $\Sigma$.
\end{definition}

It is straightforward to show that $\text{I}_{\text{aff}}$ is a well-defined quadratic form on $\Sigma$, independent of the choice of 1-adapted frame field on $\Sigma$.  As such, it may be used to define an equiaffine-invariant ``metric" on a nondegenerate surface in $\A^3$.  When $\Sigma$ is hyperbolic, equation \eqref{h-1-adapted} implies that $\text{I}_{\text{aff}}$ is equal to the indefinite quadratic form 
\[ \text{I}_{\text{aff}} = 2 \bar{\w}^1 \circ \bar{\w}^2, \] 
and so it defines a Lorentzian metric on $\Sigma$ rather than a Riemannian one.

\begin{definition}
The {\em affine Gauss curvature} $K_{\text{aff}}$ of $\Sigma$ is the Gauss curvature of the metric defined by the quadratic form $\text{I}_{\text{aff}}$.
\end{definition}

For the next step in the adaptation process, suppose that two 1-adapted frame fields $(\e_1, \e_2, \e_3), (\tilde{\e}_1, \tilde{\e}_2, \tilde{\e}_3)$ on $\Sigma$, with associated Maurer-Cartan forms $(\bar{\w}^i, \bar{\w}^i_j), (\tilde{\bar{\w}}^i, \tilde{\bar{\w}}^i_j)$, respectively, are related by a transformation of the form \eqref{1-adapted-structure-group}.  Equations \eqref{MC-forms-def} imply that 
\[
\tilde{\bar{\w}}^3_3 = \bar{\w}^3_3 + r_2 \bar{\w}^1 + r_1 \bar{\w}^2.
\]
Since $r_1, r_2$ are arbitrary real numbers, there exists a nonempty subbundle $\calF_2 \subset \calF_1$ consisting of those 1-adapted frames on $\Sigma$ for which
\begin{equation}
\bar{\w}^3_3 = 0. \label{w33-2-adapted}
\end{equation}

\begin{definition}
The bundle $\calF_2$ will be called the {\em 2-adapted frame bundle} on $\Sigma$. 
Any frame $(\e_1, \e_2, \e_3) \in \calF_2$ will be called a {\em 2-adapted frame} on $\Sigma$, and any section of $\calF_2$ will be called a {\em 2-adapted frame field} on $\Sigma$.
\end{definition}

Any two 2-adapted frames $(\e_1, \e_2, \e_3)$, $(\tilde{\e}_1, \tilde{\e}_2, \tilde{\e}_3)$ based at a point $\bx \in \Sigma$ must differ by a transformation of the form
\begin{equation}
\begin{bmatrix} \tilde{\e}_1 & \tilde{\e}_2 & \tilde{\e}_3 
\end{bmatrix} = \begin{bmatrix} \e_1 & \e_2 & \e_3 \end{bmatrix} 
\begin{bmatrix} \epsilon_1 e^{\lambda} & 0 & 0 \\[0.1in] 0 & \epsilon_2 e^{-\lambda} & 0 \\[0.1in] 0 & 0 & \epsilon_1 \epsilon_2 \end{bmatrix}, \label{2-adapted-structure-group}
\end{equation}
where $\lambda \in \R$ and $\epsilon_1, \epsilon_2 = \pm 1$.  The 2-adapted frame bundle $\calF_2$ is a principal bundle, with structure group  $G_2 \subset G_1$ equal to the subgroup of $SL(3)$ consisting of all matrices of the form in equation \eqref{2-adapted-structure-group}.

\begin{remark}
From equation \eqref{2-adapted-structure-group}, we see that the vector field $\e_3$ is now well-defined (up to sign) on $\Sigma$, independent of the choice of 2-adapted frame field on $\Sigma$.  This vector field is called the {\em affine normal vector field} on $\Sigma$. %(See \cite{NS94} for details.)
\end{remark}

Differentiating equation \eqref{w33-2-adapted} according to the structure equations \eqref{structure-eqns} implies that
\[ 0 = d\bar{\w}^3_3 = -(\bar{\w}^3_1 \& \bar{\w}^1_3 + \bar{\w}^3_2 \& \bar{\w}^2_3). \]
By Cartan's lemma, it follows that there exist functions $\ell_{ij} = \ell_{ji}$ such that
\begin{equation}
 \begin{bmatrix} \bar{\w}^1_3 \\[0.1in] \bar{\w}^2_3 \end{bmatrix} = \begin{bmatrix} 
\ell_{11} & \ell_{12} \\[0.1in] \ell_{12} & \ell_{22} \end{bmatrix} \begin{bmatrix} 
\bar{\w}^3_1 \\[0.1in] \bar{\w}^3_2 \end{bmatrix} = 
\begin{bmatrix} 
\ell_{12} & \ell_{22}  \\[0.1in] \ell_{11} & \ell_{12} \end{bmatrix} \begin{bmatrix} 
\bar{\w}^1 \\[0.1in] \bar{\w}^2 \end{bmatrix}. \label{ell-matrix}
\end{equation}

Now suppose that two 2-adapted frame fields $(\e_1, \e_2, \e_3), (\tilde{\e}_1, \tilde{\e}_2, \tilde{\e}_3)$ on $\Sigma$, with associated Maurer-Cartan forms $(\bar{\w}^i, \bar{\w}^i_j), (\tilde{\bar{\w}}^i, \tilde{\bar{\w}}^i_j)$, respectively, are related by a transformation of the form \eqref{2-adapted-structure-group}.  Equations \eqref{MC-forms-def} imply that
\[
\begin{bmatrix} \tilde{\bar{\w}}^1 \\[0.1in] \tilde{\bar{\w}}^2 \end{bmatrix} = 
\begin{bmatrix} \epsilon_1 e^{-\lambda}\, \bar{\w}^1 \\[0.1in] \epsilon_2 e^{\lambda}\, \bar{\w}^2 \end{bmatrix}, \qquad \qquad 
\begin{bmatrix} \tilde{\bar{\w}}^1_3 \\[0.1in] \tilde{\bar{\w}}^2_3 \end{bmatrix} = 
\begin{bmatrix} \epsilon_2 e^{-\lambda} \,\bar{\w}^1_3 \\[0.1in] \epsilon_1 e^{\lambda}\, \bar{\w}^2_3 \end{bmatrix} ,
\]
and it follows that the the functions $\ell_{ij}$ of equation \eqref{ell-matrix} and the corresponding functions $\tilde{\ell}_{ij}$ for the transformed frame field $(\tilde{\e}_1, \tilde{\e}_2, \tilde{\e}_3)$ are related by the equation
\begin{equation}
 \begin{bmatrix} \tilde{\ell}_{12} & \tilde{\ell}_{22} \\[0.1in] \tilde{\ell}_{11} & 
\tilde{\ell}_{12} \end{bmatrix} = \begin{bmatrix} 
\epsilon_1 \epsilon_2 \ell_{12} & e^{-2\lambda}\, \ell_{22} \\[0.1in]
e^{2\lambda}\, \ell_{11} & \epsilon_1 \epsilon_2 \ell_{12}   \end{bmatrix}   . \label{ell-action}
\end{equation}

\begin{definition}
The quadratic form
\begin{align*}
 \text{II}_{\text{aff}} & = \bar{\w}^3_1 \circ \bar{\w}^3_1 + \bar{\w}^2_3 \circ \bar{\w}^3_2 \\ & = \ell_{11} (\bar{\w}^1)^2 + 2 \ell_{12} \bar{\w}^1 \bar{\w}^2 + \ell_{22} (\bar{\w}^2)^2  
\end{align*}
on the 2-adapted frame bundle $\calF_2$ is called the {\em affine second fundamental form} of $\Sigma$.
\end{definition}

It is straightforward to show that $\text{II}_{\text{aff}}$ is a well-defined quadratic form on $\Sigma$, independent of the choice of 2-adapted frame field on $\Sigma$.

\begin{definition}
The {\em affine mean curvature} $H_{\text{aff}}$ of $\Sigma$ is defined to be $\tfrac{1}{2}$ times the trace of $\text{II}_{\text{aff}}$ with respect to the quadratic form $\text{I}_{\text{aff}}$; i.e.,  $H_{\text{aff}} = \ell_{12}.$  
\end{definition}

%As for the Euclidean analog, $H_{\text{aff}}$ is a well-defined function on $\Sigma$ up to sign.  

\begin{definition}
$\Sigma$ is called {\em affine flat} if $K_{\text{aff}}$ is identically zero on $\Sigma$, and {\em affine minimal} if $H_{\text{aff}}$ is identically zero on $\Sigma$.

\end{definition}

\begin{remark}
Unlike in Euclidean geometry, the affine Gauss curvature $K_{\text{aff}}$ is not necessarily equal to the determinant of the quadratic form $\text{II}_{\text{aff}}$. 
\end{remark}

For the remainder of this paper, we will assume that $\Sigma$ is both affine flat and affine minimal.  We will show that this assumption implies that
\[ \ell_{12} = 0, \qquad  \ell_{11} \ell_{22} = 0. \]
At each point $\bx \in \Sigma$, there are then two possibilities: either $\ell_{11} = \ell_{22} = 0$, or exactly one of $\ell_{11}, \ell_{22}$ is equal to zero.  If $\ell_{11} = \ell_{22} = 0$ at every point $\bx \in \Sigma$, then equation \eqref{MC-forms-def} implies that $d\e_3 = 0$, and hence the affine normal vector $\e_3$ is constant on $\Sigma$ and $\Sigma$ is an improper affine sphere.   On the other hand, if, say, $\ell_{22} \neq 0$ at every point of $\Sigma$, then equation \eqref{ell-action} implies that there exists a nonempty subbundle $\calF_3 \subset \calF_2$ consisting of those 2-adapted frames on $\Sigma$ for which $\ell_{22} \equiv \pm 1.$  We will not need this construction in order to obtain our normal form results in \S \ref{normal-forms-sec}, but we mention it here for the sake of completeness.

\begin{definition}
Let $\Sigma$ be a hyperbolic affine flat, affine minimal surface in $\A^3$, and suppose that $\Sigma$ contains no points where $\ell_{11} = \ell_{22} = 0$.
The bundle $\calF_3$ will be called the {\em 3-adapted frame bundle} on $\Sigma$. 
Any frame $(\e_1, \e_2, \e_3) \in \calF_3$ will be called a {\em 3-adapted frame} on $\Sigma$, and any section of $\calF_3$ will be called a {\em 3-adapted frame field} on $\Sigma$.
\end{definition}

Any two 3-adapted frames $(\e_1, \e_2, \e_3)$, $(\tilde{\e}_1, \tilde{\e}_2, \tilde{\e}_3)$ based at a point $\bx \in \Sigma$ must differ by a transformation of the form
\begin{equation}
\begin{bmatrix} \tilde{\e}_1 & \tilde{\e}_2 & \tilde{\e}_3 
\end{bmatrix} = \begin{bmatrix} \e_1 & \e_2 & \e_3 \end{bmatrix} 
\begin{bmatrix} \epsilon_1  & 0 & 0 \\[0.1in] 0 & \epsilon_2 & 0 \\[0.1in] 0 & 0 & \epsilon_1 \epsilon_2 \end{bmatrix}, \label{3-adapted-structure-group}
\end{equation}
where $\epsilon_1, \epsilon_2 = \pm 1$.  In particular, the fiber group $G_3$ of $\calF_3$ is a discrete group isomorphic to $\mathbb{Z}_2 \times \mathbb{Z}_2$, and any 3-adapted frame field on $\Sigma$ is uniquely determined by its values at any given point $\bx \in \Sigma$.

\begin{comment}

\begin{remark}
For an elliptic surface $\Sigma \subset \A^3$, the condition $K_{\text{aff}} = H_{\text{aff}} = 0$ implies that $\text{II}_{\text{aff}} = 0$, and hence  $d\e_3 = 0$ and $\Sigma$ is an improper affine sphere; in fact, it is shown in \cite{MR90} that in this case $\Sigma$ must be contained in a paraboloid.  The additional flexibility in the hyperbolic case, reflected in the possibility that $\text{II}_{\text{aff}} = \ell_{22} (\bar{\w}^2)^2 \neq 0$, is a consequence of the indefinite nature of the affine first fundamental form in the hyperbolic case.
\end{remark}

\end{comment}

\section{A local normal form}\label{normal-forms-sec}

In this section, we consider local coordinate parametrizations for $\Sigma$.  Let $U \subset \R^2$ be an open set, with coordinates $(u,v)$ on $U$, and let $\bx:U \to \A^3$ be a parametrization of a hyperbolic affine flat, affine minimal surface $\Sigma$, chosen so that the coordinate curves of $\bx$ are asymptotic curves of $\Sigma$.  (The usual Euclidean notion of an asymptotic curve for a hyperbolic surface is invariant under the group of equiaffine transformations, so this condition is well-defined.)

Define a 0-adapted frame field $(\e_1, \e_2, \e_3)$ on $\Sigma$ by setting
\[ \e_1 = \bx_u, \qquad \e_2 = \bx_v, \]
and choosing $\e_3$ to be any vector field on $\Sigma$ such that $(\e_1, \e_2, \e_3)$ is unimodular.  Then the associated Maurer-Cartan forms $(\bar{\w}^i, \bar{\w}^i_j)$ satisfy
\[ \bar{\w}^1 = du, \qquad \bar{\w}^2 = dv. \]
The condition that the coordinate curves are asymptotic is equivalent to the condition that the functions $h_{ij}$ in equation \eqref{h-matrix} satisfy 
\[ h_{11} = h_{22} = 0, \]
and therefore
\[ \bar{\w}^3_1 = h_{12}\, dv, \qquad \bar{\w}^3_2 = h_{12}\, du. \]
The condition that $\Sigma$ is nondegenerate implies that $h_{12} \neq 0$, and without loss of generality, we may assume that $h_{12}>0$: if this is not the case, simply replace $(\e_1, \e_2, \e_3)$ by the 0-adapted frame field $(\e_2, \e_1, -\e_3)$ to reverse the sign of $h_{12}$.

It is straightforward to check that the frame field
\begin{align*}
\tilde{\e}_1 & = \left(h_{12}\right)^{-(1/4)} \e_1 = \left(h_{12}\right)^{-(1/4)} \bx_u, \\
 \tilde{\e}_2 & = \left(h_{12}\right)^{-(1/4)} \e_2 = \left(h_{12}\right)^{-(1/4)} \bx_v, \\
 \tilde{\e}_3 & = \left( h_{12} \right)^{(1/2)} \e_3
\end{align*}
is 1-adapted, with Maurer-Cartan forms $(\tilde{\bar{\w}}^i, \tilde{\bar{\w}}^i_j)$ given by
\begin{alignat*}{2}
\tilde{\bar{\w}}^1 & = \left(h_{12}\right)^{(1/4)} du, & \qquad \tilde{\bar{\w}}^2 & = \left(h_{12}\right)^{(1/4)} dv, \\
\tilde{\bar{\w}}^3_1 & = \left(h_{12}\right)^{(1/4)} dv, & \qquad \tilde{\bar{\w}}^3_2 & = \left(h_{12}\right)^{(1/4)} du.
\end{alignat*}
Therefore, the affine first fundamental form \eqref{first-fundamental form} of $\Sigma$ is 
\[ \text{I}_{\text{aff}} = 2 \left(h_{12}\right)^{(1/2)} du\, dv. \]

The affine Gauss curvature of $\Sigma$ can be computed via the hyperbolic analog of Gauss's formula (see, e.g., \cite{CT80}): with $\text{I}_{\text{aff}}$ as above, we have
\[ K_{\text{aff}} = -\frac{1}{\left(h_{12}\right)^{(1/2)}} \frac{\partial^2\log \left( \left(h_{12}\right)^{(1/2)}\right)}{\partial u\partial v}. \]
The assumption that $K_{\text{aff}} = 0$ implies that
\[ \frac{\partial^2\log \left( \left(h_{12}\right)^{(1/2)}\right)}{\partial u\partial v} = 0, \]
and hence
\[ \left(h_{12}(u,v)\right)^{(1/2)} = F_1(u) F_2(v) \]
for some (nonvanishing) functions $F_1(u), F_2(v)$.  By a reparametrization of the form
\[ \tilde{u} = \int F_1(u)\, du, \qquad \tilde{v} = \int F_2(v)\, dv, \]
we can arrange that 
\[ \text{I}_{\text{aff}} = 2  d\tilde{u}\, d\tilde{v}. \]
By adjusting our frame slightly, we can construct a 1-adapted frame field $(\tilde{\e}_1, \tilde{\e}_2, \tilde{\e}_3)$ on $\Sigma$ with
\[ \tilde{\e}_1 = \bx_{\tilde{u}}, \qquad \tilde{\e}_2 = \bx_{\tilde{v}}, \]
and by adjusting $\tilde{\e}_3$, we can assume that this frame field is in fact 2-adapted.  (To reduce notational clutter, henceforth we will drop the tildes.)

The corresponding Maurer-Cartan forms $(\bar{\w}^i, \bar{\w}^i_j)$ are given by
\begin{alignat}{2}
\bar{\w}^1 & =  du, & \qquad \bar{\w}^2 & =  dv, \label{adapted-MC-forms-1} \\
\bar{\w}^3_1 & =  dv, & \qquad \bar{\w}^3_2 & =  du. \notag
\end{alignat}
Furthermore, the assumption that $H_{\text{aff}} = 0$ implies that
\begin{equation}
 \bar{\w}^1_3 = \ell_{22}\, dv, \qquad \bar{\w}^2_3 = \ell_{11}\, du. \label{adapted-MC-forms-2} 
\end{equation}
%Without loss of generality, we may assume that $\ell_{11} =0$, so that $\bar{\w}^2_3 = 0$.

In order to compute the remaining Maurer-Cartan forms, we will make use of the structure equations \eqref{structure-eqns}.  From \eqref{adapted-MC-forms-1}, we have $d\bar{\w}^1 = d\bar{\w}^2 = 0$; therefore,
\begin{align}
0 = d\bar{\w}^1 & = -(\bar{\w}^1_1 \& \bar{\w}^1 + \bar{\w}^1_2 \& \bar{\w}^2) = -(\bar{\w}^1_1 \& du + \bar{\w}^1_2 \& dv), \label{structure-comps-1} \\
0 = d\bar{\w}^2 & = -(\bar{\w}^2_1 \& \bar{\w}^1 + \bar{\w}^2_2 \& \bar{\w}^2)  = -(\bar{\w}^2_1 \& du + \bar{\w}^2_2 \& dv). \notag
\end{align}
From the relation \eqref{MC-forms-trace-relation} and the fact that $\bar{\w}^3_3=0$ for a 2-adapted frame field, it follows that $\bar{\w}^1_1 + \bar{\w}^2_2 = 0$.  Taking this into account and applying Cartan's Lemma to equations \eqref{structure-comps-1} yields
\begin{align}
\bar{\w}^1_1 & = k_1\, du + k_2\, dv, \notag \\
\bar{\w}^1_2 & = k_2\, du + k_3\, dv, \label{adapted-MC-forms-3} \\
\bar{\w}^2_1 & = k_4\, du - k_1\, dv, \notag \\
\bar{\w}^2_2 & = -k_1\, du - k_2\, dv  \notag
\end{align}
for some functions $k_1, k_2, k_3, k_4$ on $\Sigma$.

Next, from \eqref{adapted-MC-forms-1}, we have $d\bar{\w}^3_1 = d\bar{\w}^3_2 = 0$; therefore,
\begin{align}
0 = d\bar{\w}^3_1 & = -(\bar{\w}^3_1 \& \bar{\w}^1_1 + \bar{\w}^3_2 \& \bar{\w}^2_1) = 2 k_1\,du  \& dv, \label{structure-comps-2} \\
0 = d\bar{\w}^3_2 & = -(\bar{\w}^3_1 \& \bar{\w}^1_2 + \bar{\w}^3_2 \& \bar{\w}^2_2)  = 2 k_2\, du \& dv. \notag
\end{align}
Hence $k_1 = k_2 = 0$, and so $\bar{\w}^1_1 = \bar{\w}^2_2 = 0$.  Differentiating these equations yields
\begin{align}
0 = d\bar{\w}^1_1 & = -(\bar{\w}^1_2 \& \bar{\w}^2_1 + \bar{\w}^1_3 \& \bar{\w}^3_1) = k_3 k_4 \,du  \& dv, \label{structure-comps-3} \\
0 = d\bar{\w}^2_2 & = -(\bar{\w}^2_1 \& \bar{\w}^1_2 + \bar{\w}^2_3 \& \bar{\w}^3_2)  = -k_3 k_4 \, du \& dv. \notag
\end{align}
Hence $k_3 k_4 = 0$, and without loss of generality we may assume that $k_4 = 0$.  Therefore, $\bar{\w}^2_1 = 0$, and differentiating this equation yields
\[ 0 = d\bar{\w}^2_1 = -(\bar{\w}^2_1 \& \bar{\w}^1_1 + \bar{\w}^2_2 \& \bar{\w}^2_1 + \bar{\w}^2_3 \& \bar{\w}^3_1) = \ell_{11}\, du \& dv; \]
hence $\ell_{11} = 0$, and so $\bar{\w}^2_3 = 0$.  Differentiating this equation yields an identity.

Now consider the structure equation for $d\bar{\w}^1_3$: 
\[
d\bar{\w}^1_3  = -(\bar{\w}^1_1 \& \bar{\w}^1_3 + \bar{\w}^1_2 \& \bar{\w}^2_3 + \bar{\w}^1_3 \& \bar{\w}^3_3). 
\]
The left-hand side is equal to $(\ell_{22})_u\, du \& dv$, while the right-hand side is equal to zero.  Therefore, $\ell_{22} = \ell_{22}(v)$ is a function of $v$ alone.  Finally, consider the structure equation for $d\bar{\w}^1_2$:
\[
d\bar{\w}^1_2  = -(\bar{\w}^1_1 \& \bar{\w}^1_2 + \bar{\w}^1_2 \& \bar{\w}^2_2 + \bar{\w}^1_3 \& \bar{\w}^3_2). 
\]
The left-hand side is equal to $(k_3)_u\, du \& dv$, while the right-hand side is equal to $\ell_{22}\, du \& dv$.  Therefore, $(k_3)_u = \ell_{22}(v)$, and so $k_3(u,v) = u \ell_{22}(v) + f(v)$ for some function $f(v)$.

For ease of notation, let $\ell(v) = \ell_{22}(v)$.  To summarize, we have shown that the Maurer-Cartan forms associated to the 2-adapted frame $(\e_1, \e_2, \e_3)$ on $\Sigma$ are:
\begin{alignat}{3}
\bar{\w}^1 & = du, & \qquad \bar{\w}^2 & = dv, & \qquad \bar{\w}^3 & = 0, \notag \\
\bar{\w}^1_1 & = 0, & \qquad \bar{\w}^1_2 & = (u \ell(v) + f(v))\,dv, & \qquad \bar{\w}^1_3 & = \ell(v)\, dv, \label{all-MC-forms} \\
\bar{\w}^2_1 & = 0, & \qquad \bar{\w}^2_2 & = 0, & \qquad \bar{\w}^2_3 & = 0, \notag \\
\bar{\w}^3_1 & = dv, & \qquad \bar{\w}^3_2 & = du, & \qquad \bar{\w}^3_3 & = 0, \notag
\end{alignat}
where $\ell(v), f(v)$ are arbitrary functions of $v$, and that these forms satisfy the Maurer-Cartan structure equations \eqref{structure-eqns}.
%$\Sigma$ is an improper affine sphere if and only if $\ell(v) \equiv 0$.
Substituting these expressions into equations \eqref{MC-forms-def} yields the following overdetermined system of PDEs for the parametrization $\bx(u,v)$ and the 2-adapted frame field $(\e_1, \e_2, \e_3)$ on $\Sigma$:

\begin{alignat}{2}
\bx_u & = \e_1 , & \qquad \bx_v & = \e_2, \notag \\
(\e_1)_u & = 0, & \qquad (\e_1)_v & = \e_3, \label{PDE-sys} \\
(\e_2)_u & = \e_3, & \qquad (\e_2)_v & = (u \ell(v) + f(v)) \,\e_1 , \notag \\
(\e_3)_u & = 0, & \qquad (\e_3)_v & = \ell(v) \,\e_1 . \notag
\end{alignat}

The structure equations \eqref{structure-eqns} imply that the system \eqref{PDE-sys} is compatible, and the Frobenius theorem (see, e.g., \cite{IL03}) implies the following result:

\begin{theorem}
  \label{main-theorem}
  Let $U\subset \R^2$ and let $\ell(v), f(v)$ be any
  smooth, real-valued functions on $U$. Then for any
  point $(u,v)\in U$, there exists a neighborhood $V \subset U$ of
  $(u,v)$ on which the the system \eqref{PDE-sys} has a smooth solution, which defines a parametrization $\bx:V \to \A^3$ of a hyperbolic affine flat, affine minimal surface $\Sigma \subset \A^3$.  Moreover, the surface $\Sigma = \bx(V)$ is uniquely determined up to equiaffine transformations.
\end{theorem}

By an equiaffine transformation, we can assume that the functions $(\bx, \e_1, \e_2, \e_3)$ satisfy the initial conditions
\begin{equation}
 \bx(0,0) = \begin{bmatrix} 0 \\[0.05in] 0 \\[0.05in] 0 \end{bmatrix}, \qquad 
\e_1(0,0) = \begin{bmatrix} 1 \\[0.05in] 0 \\[0.05in] 0 \end{bmatrix}, \qquad 
\e_2(0,0) = \begin{bmatrix} 0 \\[0.05in] 1 \\[0.05in] 0 \end{bmatrix}, \qquad 
\e_3(0,0) = \begin{bmatrix} 0 \\[0.05in] 0 \\[0.05in] 1 \end{bmatrix}. \qquad 
 \label{ICs}
\end{equation}
Then the system \eqref{PDE-sys}, \eqref{ICs} has a unique solution in a neighborhood of $(u,v) = (0,0)$.

We can express the system \eqref{PDE-sys} explicitly as an ODE system as follows: the equations for the $u$-derivatives in \eqref{PDE-sys} imply that
\begin{align}
\bx(u,v) & = u \bar{\e}_1(v) + \bar{\bx}(v), \notag \\
\e_1(u,v) & = \bar{\e}_1(v), \label{system-soln} \\
\e_2(u,v) & = u \bar{\e}_3(v) + \bar{\e}_2(v), \notag \\
\e_3(u,v) & = \bar{\e}_3(v), \notag
\end{align}
where $\bar{\bx}(v), \bar{\e}_1(v), \bar{\e}_2(v), \bar{\e}_3(v)$ are functions of $v$ alone.  In particular, the $u$-parameter curves are straight lines in $\Sigma$, and we have the following result:

\begin{corollary}
Every hyperbolic affine flat, affine minimal surface in $\A^3$ is a ruled surface.
\end{corollary}

Substituting the expressions \eqref{system-soln} into the equations for the $v$-derivatives in \eqref{PDE-sys} yields the following ODE system for the functions $\bar{\bx}(v), \bar{\e}_1(v), \bar{\e}_2(v), \bar{\e}_3(v)$:
\begin{align}
\bar{\bx}'(v) & = \bar{\e}_2(v), \notag \\
\bar{\e}_1'(v) & = \bar{\e}_3(v), \label{ODE-sys} \\
\bar{\e}_2'(v) & = f(v)\, \bar{\e}_1(v), \notag \\
\bar{\e}_3'(v) & = \ell(v)\, \bar{\e}_1(v). \notag
\end{align}
Equations \eqref{ODE-sys} imply that $\bar{\e}_1(v)$ must satisfy the Sturm-Liouville equation
\begin{equation}
 \bar{\e}_1''(v) = \ell(v)\, \bar{\e}_1(v) \label{Sturm-Liouville-eqn}
\end{equation}
determined by the function $\ell(v)$.  Once a solution $\bar{\e}_1(v)$ to this equation has been determined, $\bar{\bx}(v)$ is obtained by integrating the equation
\begin{equation}
 \bar{\bx}''(v) = f(v)\, \bar{\e}_1(v), \label{x-eqn}
\end{equation}
taking the initial conditions \eqref{ICs} into account.

\section{Examples}\label{examples-sec}

In this section, we present some examples of hyperbolic affine flat, affine minimal surfaces; all the examples in this section are constructed by solving the system \eqref{PDE-sys} for various choices of the functions $\ell(v), f(v)$.

\begin{example}[Improper affine spheres]
If $\ell(v) \equiv 0$, then $d\e_3=0$ and $\Sigma$ is an improper affine sphere.  In this case, the system \eqref{PDE-sys}, \eqref{ICs} can be solved by quadrature, and we obtain the parametrization 
% and the 2-adapted frame $(\e_1, \e_2, \e_3)$ on $\Sigma$:
%\begin{align*}
%\bx(u,v) & = (u + F(v))\, \e_1^0 + v\, \e_2^0 + (uv + G(v))\, \e_3^0 + \bx^0,     \\
%\e_1(u,v) & =  v\,\e_3^0 + \e_1^0 ,  \\
%\e_2(u,v) & =  F'(v)\, \e_1^0 + (u + G'(v)) \e_3^0 + \e_2^0,   \\
%\e_3(u,v) & = \e_3^0, 
%\end{align*}
%where $\bx^0, \e_1^0, \e_2^0, \e_3^0$ are constant vectors with $dV(\e_1^0, \e_2^0, \e_3^0) = 1$, and the functions $F(v), G(v)$ satisfy
%\[ F''(v) = f(v), \qquad G''(v) = v f(v), \qquad F(0) = F'(0) = G(0) = G'(0) = 0.\]
%The initial conditions \eqref{ICs} imply that
%\[ \bx^0 = \begin{bmatrix} 0 \\[0.05in] 0 \\[0.05in] 0 \end{bmatrix}, \qquad 
%\e_1^0 = \begin{bmatrix} 1 \\[0.05in] 0 \\[0.05in] 0 \end{bmatrix}, \qquad 
%\e_2^0 = \begin{bmatrix} 0 \\[0.05in] 1 \\[0.05in] 0 \end{bmatrix}, \qquad 
%\e_3^0 = \begin{bmatrix} 0 \\[0.05in] 0 \\[0.05in] 1 \end{bmatrix}, \qquad 
% \]
%and then we can write
\begin{equation}
 \bx(u,v) =  \begin{bmatrix} u + F(v) \\[0.05in] v \\[0.05in] uv + G(v) \end{bmatrix}  \label{sphere-param}
\end{equation}
for $\Sigma$, where the functions $F(v), G(v)$ satisfy
\[ F''(v) = f(v), \qquad G''(v) = v f(v), \qquad F(0) = F'(0) = G(0) = G'(0) = 0.\]
Figure \ref{improper-sphere-fig-1} shows the surfaces \eqref{sphere-param} corresponding to $f(v)=0$ (the standard saddle surface $z = xy$) and $f(v)=6$.  These surfaces have parametrizations
\[ \bx(u,v) = \begin{bmatrix} u \\[0.05in] v \\[0.05in] uv \end{bmatrix}, \qquad  \bx(u,v) = \begin{bmatrix}u + 3 v^2 \\[0.05in] v \\[0.05in] uv + v^3 \end{bmatrix}, \]
respectively.
\begin{figure}[h]
\includegraphics[width=3in]{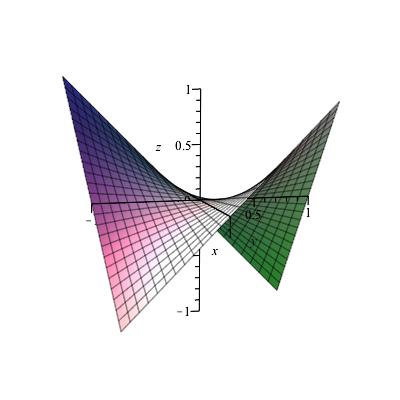}
\includegraphics[width=3in]{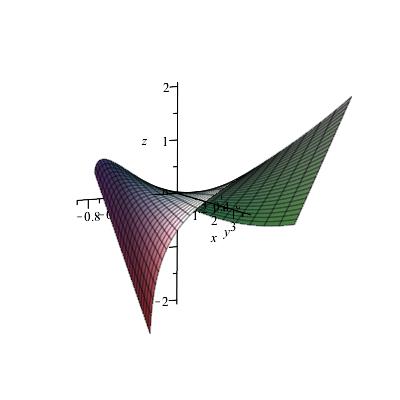}
\caption{Improper affine spheres with $f(v)=0$, $f(v) = 6$}
\label{improper-sphere-fig-1}
\end{figure}

\end{example}

\begin{remark}
The parametrization \eqref{sphere-param} describes the graph of the function
\[ z = xy + G(y) - y F(y) . \]
This agrees with the description given in \cite{MR90} of all hyperbolic, flat improper affine spheres as graphs of the form
\[ z = xy + \Phi(y), \]
where $\Phi(y)$ is an arbitrary smooth function of one variable.
\end{remark}

For the remainder of our examples, we will choose $\ell(v) \neq 0$, so that $\Sigma$ is not an improper affine sphere.

\begin{example}
Suppose that $\ell(v)$ is equal to a positive constant; i.e., $\ell(v) = a^2 > 0$.  Then the solution of equation \eqref{Sturm-Liouville-eqn} satisfying the initial conditions \eqref{ICs} is
\[ \bar{\e}_1(v) =  \begin{bmatrix} \cosh(av) \\[0.05in] 0 \\[0.05in] \frac{1}{a} \sinh(av) \end{bmatrix}, \]
and the system \eqref{PDE-sys}, \eqref{ICs} can be solved analytically to obtain the parametrization
\begin{equation}
 \bx(u,v) =  \begin{bmatrix} u \cosh(av) + F(v) \\[0.05in] v \\[0.05in] \frac{1}{a} u \sinh(av)+ G(v) \end{bmatrix}  \label{ell-pos-param}
\end{equation}
for $\Sigma$, where the functions $F(v), G(v)$ satisfy
\[
 F''(v) = f(v) \cosh(av) , \qquad G''(v) = \frac{1}{a} f(v) \sinh(av) , \qquad
 F(0) = F'(0) = G(0) = G'(0) = 0.
\]
Figure \ref{ell-pos-fig-1} shows the surfaces \eqref{ell-pos-param} corresponding to $\ell(v) = 9$ and $f(v)=0, f(v)=32 \sin(8v)$.
\begin{figure}[h]
\includegraphics[width=3in]{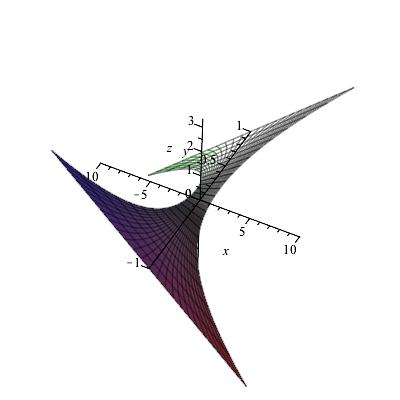}
\includegraphics[width=3in]{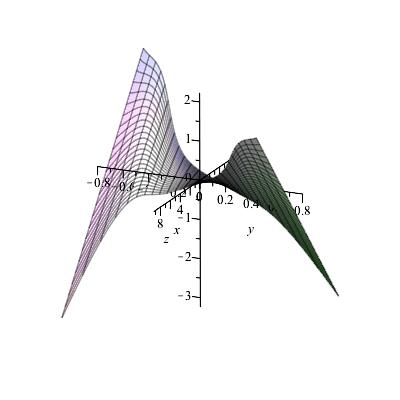}
\caption{Surfaces \eqref{ell-pos-param} with $\ell(v)=9$ and $f(v)=0$, $f(v) = 32 \sin(8v)$}
\label{ell-pos-fig-1}
\end{figure}

\end{example}

\begin{example}
Suppose that $\ell(v)$ is equal to a negative constant; i.e., $\ell(v) = -a^2 < 0$.  Then the solution of equation \eqref{Sturm-Liouville-eqn} satisfying the initial conditions \eqref{ICs} is
\[ \bar{\e}_1(v) =  \begin{bmatrix} \cos(av) \\[0.05in] 0 \\[0.05in] \frac{1}{a} \sin(av) \end{bmatrix}, \]
and the system \eqref{PDE-sys}, \eqref{ICs} can be solved analytically to obtain the parametrization
\begin{equation}
 \bx(u,v) =  \begin{bmatrix} u \cos(av) + F(v) \\[0.05in] v \\[0.05in] \frac{1}{a} u \sin(av)+ G(v) \end{bmatrix}  \label{ell-neg-param}
\end{equation}
for $\Sigma$, where the functions $F(v), G(v)$ satisfy
\[
 F''(v) = f(v) \cos(av) , \qquad G''(v) = \frac{1}{a} f(v) \sin(av) , \qquad
 F(0) = F'(0) = G(0) = G'(0) = 0.
\]
Figure \ref{ell-neg-fig-1} shows the surfaces \eqref{ell-neg-param} corresponding to $\ell(v) = -9$ and $f(v)=0, f(v)=6$.
\begin{figure}[h]
\includegraphics[width=3in]{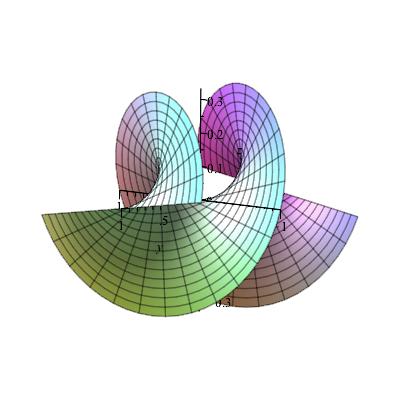}
\includegraphics[width=3in]{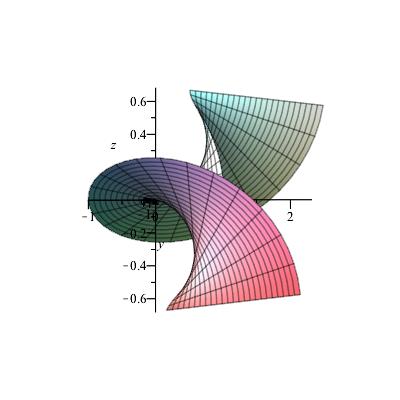}
\caption{Surfaces \eqref{ell-neg-param} with $\ell(v)=-9$ and $f(v)=0$, $f(v) = 6$}
\label{ell-neg-fig-1}
\end{figure}

\end{example}

\section{Conclusion}\label{conclusion-sec}

In affine geometry, the categories of elliptic and hyperbolic surfaces often exhibit distinctly different behavior.  As mentioned in \S \ref{intro-sec}, any elliptic affine flat, affine minimal surface in $\A^3$ must not only be an improper affine sphere, but it must in fact be contained in a paraboloid.  By contrast, there is an infinite-dimensional family of hyperbolic affine flat, affine minimal surfaces.  Magid and Ryan showed in \cite{MR90} that the improper affine spheres in this category are locally parametrized by one arbitrary function of one variable, and our results show that there is a still larger family of hyperbolic affine flat, affine minimal surfaces which are {\em not} improper affine spheres, locally parametrized by two arbitrary functions of one variable.  It would be interesting to investigate which properties of improper affine spheres may be generalized to this larger family of surfaces, and the explicit form of the PDE system \eqref{PDE-sys} should enable such investigations to be carried out fairly explicitly.

%\clearpage

\bibliographystyle{amsplain}
\bibliography{affine-flat-minimal-bib}

\end{document}